\documentclass[a4paper,12pt]{amsart}

\renewcommand{\subsection}{\section}

\usepackage{amsmath,amssymb,amscd}

\usepackage{tikz-cd}

\newtheorem{theorem}[equation]{Theorem}

\newtheorem{metaprop}[equation]{Metaproposition}

\newtheorem{proposition}[equation]{Proposition}
\theoremstyle{remark}
\newtheorem{remark}[equation]{Remark}
\theoremstyle{definition}
\newtheorem{definition}[equation]{Definition}

\newtheorem{example}[equation]{Example}

\renewcommand{\epsilon}{\varepsilon}
\renewcommand{\phi}{\varphi}

\newcommand{\cat}{\mathcal}
\newcommand{\Cat}{\ccat{Cat}}

\newcommand{\ccat}{\mathrm}

\newcommand{\widebar}{\overline}
\newcommand{\close}{\widebar}
\newcommand{\Com}{\ccat{Com}}

\newcommand{\compose}{\circ}

\newcommand{\ddeloop}{\field{\deloop}}
\newcommand{\deloop}{B}

\newcommand{\diag}{\Delta}

\renewcommand{\equiv}{\sim}

\newcommand{\equivto}{\xrightarrow{\equiv}}

\newcommand{\Family}{\ccat{Fam}}

\newcommand{\field}{\mathbb}

\newcommand{\Fun}{\ccat{Fun}}

\newcommand{\Hom}{\mathrm{Hom}}

\newcommand{\id}{\mathrm{id}}

\newcommand{\Ini}{\ccat{Init}}

\newcommand{\kore}{\textbf}

\newcommand{\Map}{\mathrm{Map}}

\DeclareMathOperator{\Ob}{Ob}

\newcommand{\op}{\mathrm{op}}

\newcommand{\Set}{\ccat{Set}}

\newcommand{\tensor}{\otimes}

\newcommand{\unity}{\boldsymbol{1}}
\newcommand{\wasreteori}{\Theta}

\title[Generalized categorification and higher theory]{A
  generalization of categorification, and higher ``theory'' of
  algebras}
\author{Matsuoka, Takuo}
\email{motogeomtop@gmail.com}



\begin{document}
\setcounter{section}{-1}

\maketitle

\begin{abstract}
We give an introduction to the topics of our forthcoming work, in
which we introduce and study new mathematical objects which we call
``higher theories'' of algebras, where inspiration for the term comes
from William Lawvere's notion of ``algebraic theory''.
Indeed, our ``theories'' are `higher order' generalizations of
coloured operad or multicategory, where we see an operad as analogous
to Lawvere's theory.

Higher theories are obtained by iterating a certain process, which we
call ``theorization'', generalizing categorification in the sense of
Louis Crane.
The hierarchy of all iterated theorizations contains in particular,
the hierarchy of all higher categories.

As an expanded introduction to the mentioned work, we here introduce
the notion of theorization, discuss basic ideas, notions, examples,
facts and problems about theorization, and describe how these lead to
our work, and what will be achieved.
\end{abstract}

\setcounter{equation}{-1}

\subsection{Higher theory}
\setcounter{subsubsection}{-1}

\subsubsection{}

This is a survey of the author's forthcoming work \cite{theory}, and
is meant to be an introduction to its topics.
In this work, we introduce and study new mathematical objects which we
call ``higher theories'' of algebras, where inspiration for the
term comes from Lawvere's notion of ``algebraic theory''
\cite{lawvere}.

Introduction of these objects leads to (among other things) a
framework for a natural explanation and a vast generalization of a
certain fact, which we would like to describe now.

A starting point for the mentioned fact is that, for a symmetric
operad $\cat{U}$ of infinity groupoids, there is a simple notion of
what we call ``$\cat{U}$-graded'' operad, which generalizes the
notions of \emph{symmetric, planar}, and \emph{braided} operad from
the cases where $\cat{U}=E_\infty, E_1, E_2$ respectively.
The notion is not complicated or hidden, but let us postpone
describing it till we are prepared shortly to reveal an
intriguing meaning of it.

The fact which can be generalized using the notion of ``higher
theory'', is the following, which is a metatheorem in the sense that
it is about a mathematical notion, rather than a mathematical object.
For an operad $\cat{U}$, let us mean by an \kore{$\cat{U}\tensor
  E_1$-monoidal} category, an associative monoidal object in the
$2$-category of $\cat{U}$-monoidal categories, or equivalently, a
$\cat{U}$-monoidal object in the $2$-category of associative monoidal
categories.

\begin{metaprop}\label{metapro:enriching-multicategory}
For every symmetric operad $\cat{U}$, the notion of $\cat{U}$-graded
operad in a symmetric monoidal category has a natural generalization
in a $\cat{U}\tensor E_1$-monoidal category.
Namely, there is a notion of $\cat{U}$-graded operad in a
$\cat{U}\tensor E_1$-monoidal category, such that the notion of
$\cat{U}$-graded operad in a symmetric monoidal category $\cat{A}$,
coincides with the notion of $\cat{U}$-graded operad in the
$\cat{U}\tensor E_1$-monoidal category underlying $\cat{A}$.
\end{metaprop}

Metaproposition naturally generalizes the familiar notions of
\begin{itemize}
\item associative algebra in a \emph{associative} monoidal category,
\item planar operad in a \emph{braided} monoidal category,
\item braided operad in a \emph{$E_3$-monoidal} infinity $1$-category
\end{itemize}
(in addition to vacuously, the notion of symmetric operad in a
symmetric monoidal category).

\begin{remark}
In order to actually have these examples, Metaproposition needs to be
interpreted in the framework of sufficiently high dimensional category
theory.
(Infinity $1$-category theory is sufficient.)
However, let us not emphasize this technical point in this survey,
even though our work will eventually be about higher category theory.
\end{remark}

\begin{remark}
The author has unfortunately failed to find a reference
for Metaproposition as stated, or any generalization of it in the
literature.
A positive is that we have found that natural ideas lead to vast
generalization of Metaproposition, even though this does not relieve
our failure of attribution.
\end{remark}

In fact, we introduce in the work much more general notion of
``grading'' (and those not just for operads), and find quite general
but natural places where the notion of graded operad (and other
things) can be enriched.
An explanation of Metaproposition
\ref{metapro:enriching-multicategory} from the general perspective to
be so acquired, will be given in Section \ref{sec:theorize-theory}.
In fact, all of these will result from extremely simple ideas, which
we would like to describe with their main consequences.

\subsubsection{}

The starting point of our work is the idea that a coloured operad or
multicategory is analogous to an algebraic theory for its role of
governing algebras over it.
In the work, we extend this by defining, for every integer $n\ge 0$,
the notion of $n$-theory,
where a $0$-theory is an algebra (commutative etc.), a $1$-theory will
be a coloured operad or ``multicategory'' (symmetric etc.), and, for
$n\ge 2$, each $n$-theory comes with a natural notion of algebra over
it, in such a way that an ($n-1$)-theory coincides
precisely with an algebra over the terminal $n$-theory, generalizing
from the case $n=1$, the fact that, e.g., a commutative algebra is an
algebra over the terminal symmeric multicategory.
Moreover, this hierarchy as $n$ varies, of ``higher'' theories,
contains the hierarchy of $n$-categories, as in fact a very small part
of it.
Quite a variety of other hierarchies of iterated categorifications,
such as operads of $n$-categories and so on, also form very small
parts of the same hierarchy.

Importance of the notion of $n$-theory lies in its significance, to be
revealed in our work, to our understanding of ($n-1$)-theories, and
hence of $0$-theories or algebras in fact, by induction.
Indeed, the notion of algebra over an $n$-theory naturally generalizes
($n-1$)-theory, and one important role of an $n$-theory is to govern
algebras over it.
In another important role, an $n$-theory provides a place in which one
can consider ($n-1$)-theories.
Namely, an $n$-theory allows one to enrich the notion of
($n-1$)-theory (including the ``generalized'' one) along it.

These two points can be combined into one expression that a
generalized (i.e., graded and enriched) ($n-1$)-theory is a
generalized (i.e., ``coloured'') functors between $n$-theories.
We shall see these points later, but in summary, the notion of
$n$-theory is important since an ($n-1$)-theory can be understood as
closely similar to a functor of $n$-theories.

In our work, we solve difficulties for precisely formulating the
notion of $n$-theory (see Section \ref{sec:theorization-general}
for more on this), and then study fundamental notions and basic facts
about these objects.
In particular, we investigate relationship among various
mathematical structures related to these objects, as well as do and
investigate various fundamental constructions.
See Sections \ref{sec:theorize-theory} and
\ref{sec:further-development} in particular.

\begin{remark}
While Lawvere's algebraic theory is about a kind of algebraic
structure which makes sense in any \emph{Cartesian} monoidal
category, an algebra over an $n$-theory makes sense (in particular) in
any, e.g., symmetric, monoidal category.
Partly for this reason, we generally call our $n$-theories higher
\emph{monoidal} theories.

It has not been clear whether there is an analogous hierarchy starting
from Lawvere's algebraic theory.
The iterated \emph{categorifications} of Lawvere's notion are
algebraic theories enriched in $n$-categories, but we are looking for
a larger hierarchy than iterated categorifications.
\end{remark}

An interesting consequence of the existence of the hierarchy of higher
theories is that, by considering an algebra over a \emph{non}-terminal
higher theory, an algebra over such a thing, and so on, we obtain
various exotic new kinds of structure, all of whom can nevertheless be
treated in a unified manner.
These structures include the hierarchy of $n$-theories for every kind
of ``grading'', specified by a choice of a (trivially graded) higher
theory $\cat{U}$ (and can in fact be exhausted essentially by all of
these).
The adjectives ``commutative'' or ``symmetric'' above refer to this
choice (the choice being the ``trivial'' grading in these cases, to be
specific).
Hierarchies in different gradings are related to each other in some
specific way which will be clarified in our work.

We expect the hierarchy would lead to new methods for studying
algebra, generalizing use of operads and multicategories, which are
just the second bit, coming next of the algebras, in this hierarchy.
In fact, the question of the existence of a similar hierarchy can be
formulated more generally, starting from much more general
kinds of ``algebraic'' structure than we have talked about so far.
Our construction of the hierarchy of $n$-theories (already of various
kinds) may be showing the meaningfulness of such a question, and this
may be our deepest contribution at the conceptual level.

The author indeed expects a similar hierarchy to exist starting from a
more general kind of algebraic structure which can be expressed
as defined by an ``associative'' operation.
Indeed, at the heart of our method is a technology of producing from a
given kind of associative operation, a new kind of associative
operation, which is based on fundamental understanding of the higher
structure of associativity.

\begin{remark}
A similar method leads to a new model \cite{category} for the theory
of higher (i.e., ``infinity infinity'') categories, including a model
of ``the infinity infinity category of infinity infinity categories''.
This will not excessively be surprising since higher theory will be a
more general kind of structure than higher category; we have already
mentioned that the hierarchy of $n$-theories contains the
hierarchy of $n$-categories.

This new model moreover has a certain convenient feature which has not
been realized on any other known model (even of infinity
\emph{$1$-categories}).
Even thought this is unfortunately not a convenient place for
describing the mentioned feature, \emph{other} features of the model
include that
\begin{itemize}
\item it is ``algebraic'' in the sense that the composition etc.~are
  given by actual operations, and
\item its construction employs only a tiny amount of combinatorics,
  and no model category theory, topology or geometry.
\end{itemize}
\end{remark}

Indeed, even though we shall discuss here only higher theoretic
structures related to algebras over multicategories, we shall also
consider in the work \cite{theory}, a modest generalization involving
higher theoretic structures related to some algebraic structures in
which the operations may have multiple inputs \emph{and} multiple
outputs, such as various versions of topological field theories.

\subsubsection{}

While the main focus of the work \cite{theory} will be on
\emph{higher} theories, Metaproposition
\ref{metapro:enriching-multicategory} can be considered as an instance
of our results in a low ``theoretic'' level of algebra.

\subsection{A generalization of categorification}
\label{sec:generalize-categorification}
\setcounter{subsubsection}{-1}

\subsubsection{}
 
The notion of $n$-theory is obtained from the notion of ($n-1$)-theory
through a certain process which generalize \emph{categorification} in
the sense of Crane \cite{crane-frenke,crane}.
In particular, the hierarchy of $n$-categories, or the iterated
categorifications of categories, is contained in the hierarchy of
$n$-theories, as has been already mentioned.
Let us start with explaining the idea for this generalized
categorification.

We would like to consider a certain process which produces a new kind
of algebraic structure from an old, so we would like to be able to
talk about \emph{kinds} of algebraic structure in general.
We would like to consider these loosely in analogy with the specific
kinds such as commutative algebra, Lie algebra, category,
$n$-category, symmetric monoidal $n$-category, symmetric operad,
planar multicategory, braided multicategory enriched in linear
$n$-categories, various versions of topological field theory, module,
and map, pairing and other combination or network of these etc.
Some of these are just algebras over multicategories, some are
categorified or similar, some are maps or functors, some make sense on
a set or on some family of sets while others require underlying linear
or other structures.

A rigourous definition of this notion will be something like Bourbaki's
\emph{species} of structure, but their particular definition is not in
a style which would be useful for us, and a rigourous definition
will not be needed here in any case.
(In Section \ref{sec:theorization-general}, we shall formalize some
feature of the notion just in order to help ourselves to consider some
further notions, but without any attempt to be comprehensive.)

\subsubsection{}

Now, for a kind of algebraic structure definable on a set, its
\emph{categorification} is an analogous algebraic structure which
makes sense on a category, but is defined replacing structure maps by
functors, and structural equations by suitably coherent isomorphisms,
forming a part of the structure.
A basic feature expected of the categorified structure is that, if a
category $\cat{C}$ is equipped with a categorified form of a certain
kind of algebraic structure, then the original, uncategorified form of
the same structure should naturally make sense on any object of
$\cat{C}$.
Indeed, a commutative algebra makes sense in a symmetric monoidal
category, and an associative algebra makes sense in any
associative monoidal category.

However, this is not the most general instance of the phenomenon.
In a symmetric monoidal category for example, the notion of algebra
makes sense over any symmetric operad or multicategory, and the same
moreover makes sense also \emph{in} any symmetric multicategory.
Indeed, an algebra over a symmetric multicategory $\cat{U}$ in a
symmetric multicategory $\cat{V}$ is simply a morphism
$\cat{U}\to\cat{V}$.
Similarly, the notion of associative algebra, as well as the notion
of algebra over any planar multicategory, makes sense in any
associative monoidal category, and more generally, in any planar
multicategory in the same manner.

While the notion of associative monoidal category categorifies the
notion of associative algebra, there is a more general process than
that of categorification which produces the notion of planar
multicategory from the notion of associative algebra, and symmetric
multicategory from commutative algebra.
In general, this process produces from a given kind of algebraic
structure, a new kind of algebraic structure generalizing its
categorification, in such a manner that the original notion of algebra
reduces to the notion of algebra over the terminal one among the new
objects defined (meaning symmetric multicategories, for commutative
algebras, so generalizing the simple fact that an commutative algebra
is an algebra over the terminal symmetric multicategory).

Let us thus recall how one may naturally arrive at the notion of
symmetric operad starting from the notion of commutative algebra.
(And we are suggesting that the same procedure will produce the notion
of planar operad out of the notion of associative algebra, for
example.)
Specifically, let us try to find the notion of symmetric operad out of
the desire of generalizing the notion of commutative algebra to the
notions of certain other kinds of algebra which makes sense in a
symmetric monoidal category.
Indeed, one of the most important role of a multicategory is
definitely the role of governing algebras over it.

The way how we generalized the notion of commutative algebra is as
follows.
Namely, the structure of a commutative algebra on an object $A$ of a
symmetric monoidal category, can be given by a single $S$-ary
operation $A^{\tensor S}\to A$ for every finite set $S$ which,
collected over all $S$, has appropriate consistency.
We get a generalization of this by allowing not just a single $S$-ary
operation, but a \emph{family} of $S$-ary operations parametrized by a
set or a space prescribed for $S$.
(One may use, instead of spaces, any suitable mathematical objects
here depending on the category in which the object $A$ lives.
In fact, results of our work will enable us to recognize a wide
variety of mathematical objects which one can use for the `enrichment'
of the structure here, or in much more various situations.)
This ``space of $S$-ary operations'' for each $S$, is the first bit of
the data defining an operad.
Having this, we next would like to compose these operations just as we
can compose multiplication operations of a commutative algebra,
and the composition should have appropriate consistency.
A symmetric multicategory is simply a more general version of this,
with many objects, or ``colours''.

A similar procedure can be imagined once a kind of ``algebraic''
structure in a very broad sense, in place of commutative or
associative algebra, is specified as a specific kind of system of
operations.
Inspired by Lawvere's notion of an \emph{algebraic theory}
\cite{lawvere}, we call a multicategory, emphasizing its role of
governing a particular kind of algebras, also a (symmetric)
\kore{$1$-theory}, and then generally call \kore{theorization}, a
process similar to the process above through which we have obtained
$1$-theories from the notion of commutative algebra.
The result of such a process will also be called a \kore{theorization}.
Thus, the notion of $1$-theory is a \emph{theorization} of the notion
of commutative algebra.
We shall see in Sections \ref{sec:theorization-general} and
\ref{sec:basic-construction}, how the process of theorization indeed
generalizes the process of categorification.

\subsection{Theorization of algebra}
\label{sec:theorize-algebra}
\setcounter{subsubsection}{-1}

\subsubsection{}

As the simplest example of a theorization process next to the one
which we have seen in the previous section, let us consider
theorization of the
notion of $\cat{U}$-algebra for a symmetric operad $\cat{U}$.
By using the same method as in the previous section, we shall obtain a
theorization of the notion of $\cat{U}$-algebra, which we call
\emph{$\cat{U}$-graded operad}.
Let us assume for simplicity, that $\cat{U}$ is an uncoloured operad.

Recall that the structure of a $\cat{U}$-algebra $A$ is determined by
an action of every operator $u$ belonging to $\cat{U}$ on $A$.
If $u$ belongs to the set ``of $S$-ary operations'' in $\cat{U}$ for a
finite set $S$, then it should act as an $S$-ary operation $A^{\tensor
  S}\to A$.
Now, to theorize the notion of $\cat{U}$-algebra given by an action of
the operators of $\cat{U}$, means to modify the definition
of this structure by replacing an action of every operator $u$ in
$\cat{U}$, by a choice of the set ``of operations of shape (so to
speak) $u$''.
We call an element of this set an \kore{operation of degree $u$}.

Thus the data of a $\cat{U}$-graded operad $\cat{X}$ should include,
for every operation $u$ in $\cat{U}$, a set whose element we shall
call an operation in $\cat{X}$ of degree $u$.
If the operation $u$ is $S$-ary in $\cat{U}$, then we shall say that
any operation of degree $u$ in $\cat{X}$ has \kore{arity} $S$.

There should further be given a consistent way to compose the
operations in $\cat{X}$ which moreover respects the degrees of the
operations.
These will be a complete set of data for a \kore{$\cat{U}$-graded
  operad} $\cat{X}$.

There is actually a coloured version of this which we call
$\cat{U}$-graded \kore{multicategory} or \kore{$1$-theory}, so this is
a more general theorizations of $\cat{U}$-algebra.
$\cat{U}$-graded multicategory is in fact also a generalization of
$\cat{U}$-monoidal category, generalizing the fact that symmetric
multicategory was a generalization of symmetric monoidal category.

\subsubsection{}

By reflecting on what we have done above, we immediately find that a
$\cat{U}$-graded operad is in fact exactly a symmetric operad
$\cat{Y}$ equipped with a morphism $P\colon\cat{Y}\to\cat{U}$.
The relation between $\cat{X}$ above and $\cat{Y}$ here is that an
$S$-ary operation in $\cat{Y}$ is an $S$-ary operation in $\cat{X}$
of \emph{arbitrary} degree.
The map $P$ maps an operation in $\cat{Y}$ to the degree it had in
$\cat{X}$.
Conversely, given an $S$-ary operation $u$ in $\cat{U}$, an operation
in $\cat{X}$ of degree $u$ is an $S$-ary operation in $\cat{Y}$ which
lies over $u$.
In particular, $\cat{U}$ can be identified with the terminal
$\cat{U}$-graded operad, which has exacly one operation of each
degree.

A $\cat{U}$-graded multicategory also turns out to be just a symmetric
multicategory equipped with a map to $\cat{U}$.
Now, given a $\cat{U}$-graded operad $\cat{X}$, an $\cat{X}$-algebra
in a $\cat{U}$-graded $1$-theory $\cat{Y}$ will be just a functor
$\cat{X}\to\cat{Y}$ of $\cat{U}$-graded $1$-theories.

\begin{example}\label{ex:initially-graded}
A multicategory graded by the initial operad $\Ini$ is a multicategory
with only unary multimaps, which is equivalent to a category.
\end{example}

A similar theorization of the notion of $\cat{U}$-algebra, can also be
defined for a \emph{coloured} symmetric operad $\cat{U}$, and a
$\cat{U}$-graded $1$-theory will be again a symmetric multicategory
equipped with a functor to $\cat{U}$.

\begin{example}\label{ex:graded-by-category}
Recall, as noted in Example \ref{ex:initially-graded}, that a category
$\cat{C}$ can be considered as a symmetric
multicategory having only unary maps.
If we consider $\cat{C}$ as a multicategory in this way, then a
$\cat{C}$-graded $1$-theory is a category equipped with a functor to
$\cat{C}$, and this theorizes $\cat{C}$-algebra, or functor on
$\cat{C}$ (which one might also call a left $\cat{C}$-module).
On the other hand, a categorification of $\cat{C}$-algebra (in the
category of sets) is a category valued functor $\cat{C}\to\Cat$, and
among the theorizations, categorifications correspond to op-fibrations
over $\cat{C}$.
\end{example}

\begin{remark}\label{rem:groupoid-terminology}
This is a technical remark.

In this example (and everywhere else in this survey), the category of
categories should be considered as enriched over groupoids.
See Section \ref{sec:theorization-general}.
Then a functor $\cat{C}\to\Cat$, should be understood as a functor in
the usual ``weakened'' sense (which is sometimes called a
\emph{pseudo}-functor).

All other categorical terms in this survey should be understood in the
similar manner when there is enrichment of the relevant categorical
structures in groupoids.
(The reader who is comfortable with homotopy theory may instead
replace all sets/groupoids with infinity groupoids, and understand
everything as enriched in infinity groupoids.)
\end{remark}

Suppose given a category $\cat{C}$ and two functors
$F,G\colon\cat{C}\to\Cat$, corresponding respectively to categories
$\cat{X},\cat{Y}$ lying over $\cat{C}$, mapping down to $\cat{C}$ by
op-fibrations.
Note that, by Example \ref{ex:graded-by-category}, $F$ and $G$ are
categorified $\cat{C}$-modules, and $\cat{X}$ and $\cat{Y}$ as
categories over
$\cat{C}$, are the corresponding $\cat{C}$-graded $1$-theories.

In this situation, the relation between maps $F\to G$ and maps
$\cat{X}\to\cat{Y}$, is as follows.
Namely, a functor $\phi\colon\cat{X}\to\cat{Y}$ of categories over
$\cat{C}$ (see Remark \ref{rem:groupoid-terminology}, to be
technical), corresponds to a map $F\to G$ if and only if $\phi$
preserves coCartesian maps, and an arbitrary functor $\phi$ over
$\cat{C}$ only corresponds to a \emph{lax} map $F\to G$. (defined with
the \emph{$2$-category} structure of $\Cat$ in consideration).

A similar pattern can be observed on theorization in general.

\subsection{Theorization in general}
\label{sec:theorization-general}
\setcounter{subsubsection}{-1}

\subsubsection{}

For the idea for theorization of a more general algebraic structure,
the notion of \emph{profunctor}/\emph{distributor}/\emph{bimodule}
is useful.
For categories $\cat{C},\cat{D}$, a
\kore{$\cat{D}$--$\cat{C}$-bimodule} is a functor
$\cat{C}^\op\times\cat{D}\to\Set$, where $\Set$ denotes the
category of sets.
The category of $\cat{D}$--$\cat{C}$-bimodules contains the opposite
of the functor category $\Fun(\cat{C},\cat{D})$
as a full subcategory, where a functor $F\colon\cat{C}\to\cat{D}$ is
identified with the bimodule $\Map_\cat{D}(F-,-)$.
Let us say that this bimodule is \kore{corepresented} by $F$.
By symmetry, the category of bimodules also contains
$\Fun(\cat{D},\cat{C})$.
However, for the purpose of theorization, we treat $\cat{C}$ and
$\cat{D}$ asymmetrically, and mostly consider only
\emph{co}representation of bimodules.
Bimodules compose by tensor product, to make categories form a
$2$-category, extending the $2$-category formed with (opposite)
functors as $1$-morphisms, by the identification of a functor with the
bimodule corepresented by it.

\subsubsection{}

Now, we would like to define theorization of algebraic structures
definable on a set, so the notion of $\cat{U}$-graded multicategory we
have arrived at in the previous section will indeed be a theorization
in this sense, of $\cat{U}$-algebra in the category $\Set$, or
``$\cat{U}$-\emph{monoid}''.
Before beginning, let us mention that our purpose for giving a
definition is merely to make the idea visible, and is \emph{not} to
analyse the notion at an abstract level.
We just would like to point out a pattern about some concrete
examples which interests us.
Therefore, it will not be a precisely adjusted and most general
definition which will be useful for us, or will be desired for.

Still, a workable definition would need \emph{some} formulation of
the notion of structure.
Let us define as follows, just for the discussions here.

\begin{definition}\label{def:kind-of-structure}
Let $\cat{C}$ be a category.
Then a \kore{kind} of structure on an object of $\cat{C}$, is given by
a pair consisting of
\begin{itemize}
\item a category $\cat{K}$, to be called the category ``of
  \kore{structured} objects'', and
\item a functor $F\colon\cat{K}\to\cat{C}$, to be called the
  \kore{forgetful} functor.
\end{itemize}
We call a structure of the kind specified by these data a
\kore{$\cat{K}$-structure} (while remembering $F$).
\end{definition}

For example, we would like to consider structured categories along
this idea.
However, there is one thing to bear in mind in this case.
Namely, when we apply this idea to $\cat{C}=\Cat$, we would like to
consider the category $\Cat$ of categories (with a fix limit for size)
as enriched in groupoids.
Namely, for categories $X,Y\in\Cat$, we consider $\Map_\Cat(X,Y)$ as
the groupoid formed by functors $X\to Y$ and isomorphisms between
those functors.
In Definition \ref{def:kind-of-structure}, if the category $\cat{C}$
is enriched in groupoids (like $\Cat$), then we normally consider
$\cat{K}$ and $F$ which are also enriched in groupoids (see Remark
\ref{rem:groupoid-terminology}, to be technical).

In some other cases, we would like to consider also structures
definable on a family indexed by some collection, of sets.
For example, for a category $X$, there is the family
$\Map_X:=\big(\Map_X(x,y)\big)_{x,y}$ of sets of morphisms,
parametrized by the pairs $x,y$ of objects of $X$, so the structure of
$X$ can be understood as defined on this family $\Map_X$ of sets, by
the composition operations.
In general, when the indexing collection $\Lambda$ of a family is
chosen and fixed for our consideration, let us, for a category
$\cat{C}$ (possibly enriched in groupoids), denote by
$\Family_\Lambda(\cat{C})$, the (enriched) category of families of
objects of $\cat{C}$ indexed by $\Lambda$.

Now recall that a set can be identified with a homotopy $0$-type,
namely, a groupoid in which every pair of maps $f,g\colon x\equivto y$
between the same pair of objects, are equal.
With this identification, every category can be considered as enriched
in groupoids since it is enriched in homotopy $0$-types.
Using this enrichment, the category $\Set$ of sets can be considered
as embedded into the category $\Cat$ enriched in groupoids, as the
full subcategory consisting of homotopy $0$-types.

Let us now define as follows, \emph{without} excluding trivial
examples.

\begin{definition}\label{def:categorification-example}
Suppose given a kind of structure $\cat{K}$ definable on a family
indexed by a collection $\Lambda$, of sets, with forgetful functor
$F\colon\cat{K}\to\Family_\Lambda(\Set)$.
Then a kind of structure $\close{\cat{K}}$ on family of categories
indexed by $\Lambda$, with forgetful functor
$\close{F}\colon\close{\cat{K}}\to\Family(\Cat)$ (of categories
enriched in groupoids) is a \kore{categorification} of
$\cat{K}$-structure, if there is given a Cartesian square
\[\begin{tikzcd}
\cat{K}\arrow{d}[swap]{F}\arrow{r}
&\close{\cat{K}}\arrow{d}{\close{F}}\\
\Family_\Lambda(\Set)\arrow[hook]{r}&\Family_\Lambda(\Cat).
\end{tikzcd}\]
Namely, we say that such a square exhibits $\close{\cat{K}}$-structure
\emph{as} a categorification of $\cat{K}$-structure.
\end{definition}

In other words, a \emph{categorification} of $\cat{K}$-structure is a
generalization of $\cat{K}$-structure to structure on a family of
categories.

\begin{remark}
The idea just stated is actually more important than Defintion.
In fact, Definition should rather be seen as an example of a situation
where the stated idea can be applied, since examples we encounter in
practice can often be ``categorifications'' in suitably more general
sense.
However, what would be important for our purposes will be to
understand concrete situations based on the stated idea, and the
author does not expect to gain much, as far as our work is concerned,
by formulating the notion more generally.
The situation covered by our definition can never the less be
considered typical, so we shall be content with it here.

A similar remark will apply also to our definitions of related notions
to follow below in this section.

Note also that niether the idea above nor Definition excludes trivial
or unnatural examples.
We shall not mind this since we do have interesting examples such as
the hierarchy of $n$-categories.
\end{remark}

\begin{remark}
In practice, the context is often not categorical, but \emph{infinity}
$1$-categorical.
In such a context, one usually consider the infinity $1$-category of
groupoids in place of the category of $\Set$.
In this situation, a ``categorification'' of a kind of structure on
infinity groupoid means a kind of structure generalizing it on
infinity $1$-category.
\end{remark}

\subsubsection{}

Next, we would like to suppose given a kind of structure on family of
categories, indexed by a specific collection, say $\Lambda$.
For example, we can consider the structure of a lax associative
monoidal category as a structure on its underlying category, say $X$.
In this case, the indexing collection $\Lambda$ is a set consisting of
one point, and the structure is given, for example, by the
``multiplication'' functor $\tensor\colon X^n\to X$ for every integer
$n\ge 0$, and associativity maps between suitable combinations of
these functors, which satisfy suitable equations.
For another example, if we would like to consider instead, an
associative monoidal structure which is not lax, then we would leave
$\Lambda$ the same, and the structure functors $\tensor$ similar, but
would require the associativity maps to be isomorphisms.

Similarly, if we consider any specific kind of structure on family of
categories, then at least part of the structure may be defined by
\begin{itemize}
\item functors or bimodules/distributors/profunctors (like the
  monoidal multiplication functors above) between some specific pairs
  of categories in the family, or of the direct products of some
  specific members of the family, as well as
\item maps (like the associativity maps above) between
  functors/bimodules obtained by combining some of those structure
  functors/bimodules in a specific manner (which we shall call
  \kore{structure $2$-maps}),
  or similar isomorphism, which are required to satisfy
\item some particular equations.
\end{itemize}

\begin{definition}\label{def:relax-virtualize}
Fix a collection $\Lambda$ for indexing families, and consider the
groupoid-enriched category $\Family(\Cat)=\Family_\Lambda(\Cat)$
formed by families of categories indexed by $\Lambda$.
Suppose given two kinds $\cat{K}$, $\close{\cat{K}}$ of structure on
object of $\Family(\Cat)$, with respective forgetful
(enriched) functors $F\colon\cat{K}\to\Family(\Cat)$ and
$\close{F}\colon\close{\cat{K}}\to\Family(\Cat)$.
Suppose further given a functor $T\colon\cat{K}\to\close{\cat{K}}$ of
enriched categories over $\Family(\Cat)$ (see Remark
\ref{rem:groupoid-terminology}, to be technical).
Then
\begin{itemize}
\item $\close{\cat{K}}$-structure is said to be a \kore{relaxation} of
  $\cat{K}$-structure if $T$ is fully faithful on the underlying
  $2$-groupoids, and $\cat{K}$-structures are characterized among
  $\close{\cat{K}}$-structures as those in which some specific
  structure $2$-maps (see above) are invertible,
\item $\close{\cat{K}}$-structure is said to be a
  \kore{virtualization} of $\cat{K}$-structure if $T$ is fully
  faithful on the underlying $2$-groupoids, and $\cat{K}$-structures
  are characterized among $\close{\cat{K}}$-structures as those in
  which some specific structure bimodules are corepresentable.
\end{itemize}
Namely, we say that $T$ exhibits $\close{\cat{K}}$-structure \emph{as}
a relaxation or virtualization if one of these conditions are satisfied.
\end{definition}

In other words, a \emph{relaxation} is a generalization of the
structure by allowing non-invertible maps in place of some specific
structure \emph{isomorphisms}, and a \emph{virtualization} is a
generalization by allowing non-corepresentable bimodules in place of
some specific structure \emph{functors}.

\subsubsection{}

The idea of theorization we have described in the previous sections,
can now be expressed as follows.

\begin{definition}\label{def:theorization-example}
Suppose given a kind $\cat{K}$ of structure on family of sets indexed
by a collection $\Lambda$.
Then a \kore{theorization} of $\cat{K}$-structure is a virtualization
of a relaxation of a categorification of $\cat{K}$-structure.
\end{definition}

It is obvious from the definitions that none of categorification,
relaxation or virtualization is uniquely determined from the original
kind of structure, or is guaranteed to exist in a non-trivial manner.
Therefore, given a specific kind of structure, it seems to be a
usually non-trivial question whether the kind of structure has any
good, or even just non-trivial, theorization.

\begin{remark}\label{rem:colour-after-theorization}
The notion of theorization of Definition is the ``coloured'' version
which we did not discuss in Section
\ref{sec:generalize-categorification} or \ref{sec:theorize-algebra}.
A \emph{colour} in the theorized structure is an object of a category
in the family indexed by $\Lambda$.
This generalizes the colours in a multicategory.
See Example \ref{ex:graded-as-theorization} below.
\end{remark}

\begin{example}\label{ex:graded-as-theorization}
For a multicategory $\cat{U}$, let us try to interpret
$\cat{U}$-graded multicategory as a ``theorization'' in this sense.

Firstly, we consider the structure of a $\cat{U}$-monoid as a
structure on family of sets indexed by the objects of $\cat{U}$.
Namely, we consider a $\cat{U}$-monoid $A$ as consisting of
\begin{itemize}
\item for every object $u\in\cat{U}$, a set $A(u)$,
\item for every finite set $S$ and an $S$-ary operation $v\colon
  u\to u'$ in $\cat{U}$, where $u=(u_s)_{s\in S}$ is a family of
  objects of $\cat{U}$ indexed by $S$, a map $A(v)\colon A(u)\to
  A(u')$, where $A(u):=\prod_{s\in S}A(u_s)$,
\end{itemize}
and then consider the latter as a structure on the family
$\Ob A:=\big(A(u)\big)_{u\in\cat{U}}$ of sets indexed by the objects
$u$ of $\cat{U}$.
Namely, the structure is the action of every multimaps $v$ in
$\cat{U}$ on the relevant members of the family $\Ob A$.

In order to indeed obtain a similar family \emph{of categories} from a
$\cat{U}$-graded multicategory $\cat{X}$, we take for every object $u$
of $\cat{U}$, the category formed by the objects of $\cat{X}$ of
degree $u$, and maps (i.e., unary multimaps) between them of degree
$\id_u$.
Namely, if we denote this category by $\cat{X}_u$, then we think of
the family $\Ob\cat{X}:=(\cat{X}_u)_u$ of categories indexed by
objects $u\in\cat{U}$, as underlying $\cat{X}$.

The rest of the structure of $\cat{X}$ can then be considered as the
lax associative action of the rest of the multimaps $f$ in $\cat{U}$,
on these categories $\cat{X}_u$, each $f$ acting as the bimodule
formed by the multimaps in $\cat{X}$ of degree $f$, so $\cat{X}$ can
be interpreted as obtained by putting a theorized $\cat{U}$-algebra
structure on the family $\Ob\cat{X}$.
\end{example}

If a $\cat{U}$-graded multicategory $\cat{X}$ is seen as a theorized
structure in this manner, then a colour in this theorized structure is
an object of a category $\cat{X}_u$, where $u$ is any object of
$\cat{U}$.
In other words, it is an object of the multicategory $\cat{X}$.

\begin{remark}\label{rem:action-of-identity}
This will be minor and technical.

However, a natural theorization which Definition expects of
$\cat{U}$-monoid would appear to be lax $\cat{U}$-algebra in the
$2$-category mentioned above formed by categories and bimodules
between them.
This does not coincide with our desired theorization, which is
$\cat{U}$-graded multicategory.

Indeed, for an object $u$ of $\cat{U}$, if a category, say
$\cat{X}_u$, is associated to $u$, and the identity map of $u$ acts on
$\cat{X}_u$ in our $2$-category of bimodules, then this acton gives
another category, say $\cat{Y}_u$, with objects the objects of
$\cat{X}_u$, and a map, say $F\colon\cat{X}_u\to\cat{Y}_u$, of the
structures of categories on the same collection of objects.
However, the theorization in the idea described in the previous
sections, is not where $\id_u$ acts on an already existing category
$\cat{X}_u$, but where the structure of $\cat{X}_u$ itself as a
category, is the action of $\id_u$.

In other words, we usually do not just want to consider a relaxed
structure in the $2$-category of categories and bimodules, but we
would further like to require that the resulting map corresponding to
$F$ in the example above, associated to each of the
categories in the family, to be an isomorphism.
\end{remark}

\begin{remark}
In Example \ref{ex:graded-as-theorization}, there is another category
structure on the objects of $\cat{X}_u$, in which a map is a (unary)
map in $\cat{X}$ of degree an \emph{arbitrary} endomorphism of $u$ in
$\cat{U}$, rather than just the identity.
This of course does not interfere with Remark
\ref{rem:action-of-identity}.
\end{remark}

\subsection{A basic construction}
\label{sec:basic-construction}
\setcounter{subsubsection}{-1}

Definitions \ref{def:theorization-example} and
\ref{def:relax-virtualize} imply that a categorified structure is an
instance of the theorized form of the same structure.
Given a categorified structure $\cat{X}$, let us denote by
$\wasreteori\cat{X}$, the theorized structure corresponding to
$\cat{X}$.
Concretely, $\wasreteori\cat{X}$ is obtained by replacing as needed,
structure functors of $\cat{X}$ with bimodules copresented by them.
The construction $\wasreteori$ generalizes the usual way to construct
a multicategory from a monoidal category.
Let us say that $\wasreteori\cat{X}$ is \kore{represented} by
$\cat{X}$.

\begin{remark}
Our definition \ref{def:theorization-example} of theorization allowed
theorized structures to have colours (Remark
\ref{rem:colour-after-theorization}).
This flexibility is playing an important role here.
Namely, if the categorified structure $\cat{X}$ is structured on a
family, say $\Ob\cat{X}=(\Ob_\lambda\cat{X})_{\lambda\in\Lambda}$ of
categories, where $\Lambda$ denotes the collection indexing our
families, then any object of $\Ob_\lambda\cat{X}$ for
$\lambda\in\Lambda$ is being a
colour in the theorized structure $\wasreteori\cat{X}$.

Thus the coloured version of the notion of theorization is indeed
necessary in order for every
categorified structure to be included in theorized structures.
\end{remark}

As in the case of monoidal structure, $\wasreteori$ is usually only
faithful, but not full.
Indeed, for (families of) categories $\cat{X},\cat{Y}$ equipped with
categorified structures, a morphism
$\wasreteori\cat{X}\to\wasreteori\cat{Y}$ is equivalent to a
\emph{lax} morphism $\cat{X}\to\cat{Y}$.

\begin{example}
Let $\cat{U}$ be a symmetric multicategory, and let $\cat{A}$ be a
$\cat{U}$-monoidal category.
Then a $\cat{U}$-algebra in $\cat{A}$ is the same thing as a
$\cat{U}$-algebra in the $\cat{U}$-graded $1$-theory
$\wasreteori\cat{A}$, which is equivalent to a map to
$\wasreteori\cat{A}$ from the terminal $\cat{U}$-graded $1$-theory.

However, the terminal $\cat{U}$-graded $1$-theory is represented by
the unit $\cat{U}$-monoidal category $\unity$, so this is in agreement
with the familiar fact that a $\cat{U}$-algebra in $\cat{A}$ is the
same thing as a lax $\cat{U}$-monoidal functor $\unity\to\cat{A}$.
\end{example}

\begin{remark}
The functor $\wasreteori$ has a left adjoint (which in fact can be
described in a very concrete manner).
In the example above, $\cat{U}$-algebra in $\cat{A}$ is thus
equivalent to a $\cat{U}$-monoidal functor to $\cat{A}$ from the
$\cat{U}$-monoidal category freely generated from the terminal
$\cat{U}$-graded $1$-theory $\wasreteori\unity$ (which thus has a
concrete description).
\end{remark}

\subsection{Theorization of category}
\label{sec:theorize-category}
\setcounter{subsubsection}{-1}

\subsubsection{}

For illustration of the general definition, let us describe a natural
theorization of the notion of category, which we shall call
``categorical theory'' here.

We have mentioned in Section \ref{sec:theorization-general}, that a
category may be considered as an algebraic structure on a family of
sets ``of morphisms'', parametrized by pairs of its objects, and, if
we choose and fix a collection as the collection of objects for our
categories, then $2$-category with the same collection of objects, is
a ``categorification'' of those categories in the sense of Definition
\ref{def:categorification-example}.
Categorical theory will be a theorization of category in the sense of
Definition \ref{def:theorization-example}, whose associated
categorification is $2$-category.

The description of a \emph{categorical theory} is as follows.
Firstly, it, like a $2$-category (our categorification) has objects,
$1$-morphisms, and sets of $2$-morphisms.
$1$-morphisms do not compose, however.
Instead, for every nerve
$f\colon x_0\xrightarrow{f_1}\cdots\xrightarrow{f_n}x_n$ of
$1$-morphisms and a $1$-morphism $g\colon x_0\to x_n$, one has the
notion of \emph{$n$-ary $2$-multimap} $f\to g$.
The $2$-morphisms, which were already mentioned, are just unary
$2$-multimaps.
There are given unit $2$-morphisms and associative composition for
$2$-multimaps, analogously to the similar operations for multimaps in
a planar multicategory.

A $2$-category is in particular a categorical theory, in which a
$2$-multimap $f\to g$ is a $2$-morphism $f_n\compose\cdots\compose
f_1\to g$ in the $2$-category.
Between two $2$-categories, a natural map \emph{of categorical
  theories} is not precisely a functor, but is a lax functor of the
$2$-categories.

As we have also suggested, a categorical theory is also a
generalization of a planar multicategory.
Indeed, planar multicategory was a theorization of associative monoid.
The relation between the notions of planar multicategory and of
categorical theory, is parallel to the relation between the notions of
associative monoid and of category.
Namely, categorical theory is a `many objects' (or ``coloured'')
version of a planar multicategory, where the word ``object'' here
refers to one at a deeper level than the many objects which
a multicategory (as a ``coloured'' operad) may already have are at.

One thing one should note then, is that, while the $2$-multimaps
in a categorical theory is generalizing the multimaps in a planar
multicategory here, the $1$-morphisms in a categorical theory is
generalizing the \emph{objects} of a planar multicategory, and no
longer have the characteristic of operators like the $1$-morphisms in
a category.
Indeed, a $1$-morphism in a categorical theory and an object of a
planar multicategory are both ``colours'' in the sense of Remark
\ref{rem:colour-after-theorization}, and we have also mentioned
earlier that there is no operations of composition given for
$1$-morphisms in a categorical theory.

What we said above in comparison of the structures of a categorical
theory and of a planar multicategory, is that a categorical theory
$\cat{C}$ has one more layer of `colouring' under the $1$-morphisms,
given by the collection of the objects of $\cat{C}$.

\subsubsection{}

Following the general pattern about theorization, there is a notion of
category in a categorical theory, and, as an uncoloured version of it,
monoid in a categorical theory, which generalizes a monad in a
$2$-category.
A monad in a $2$-cateogry $\cat{C}$ was a lax functor to $\cat{C}$
from the terminal $2$-category, which can also be considered as a map
between the categorical theories represented by these $2$-categories.
See Section \ref{sec:basic-construction}. 
However, the latter is a monoid in the target categorical theory
$\wasreteori\cat{C}$ by definition.
More generally, a category in a categorical theory can be described as
a coloured version of a map of categorical theories.

\begin{example}\label{ex:theory-from-monad}
Let $\cat{C}$ be a category enriched in groupoids, and let $M$ be a
monad \emph{on} $\cat{C}$.
Then there is a categorical theory as follows.
\begin{itemize}
\item An object is an object of $\cat{C}$.
\item A map $x\to y$ is a map $Mx\to y$ in $\cat{C}$.
\item Given a sequence of objects $x_0,\ldots,x_n$ and maps $f_i\colon
  Mx_{i-1}\to x_i$ in $\cat{C}$ and $g\colon Mx_0\to x_n$, the set
  $\Hom(f,g)$ of $2$-multimaps $f\to g$, is the set of
  \emph{commutative} diagrams
  \[\begin{tikzcd}[column sep=large]
    M^nx_0\arrow{d}[swap]{m}\arrow{r}{M^{n-1}f_1}
    &\cdots\arrow{r}{f_n}
    &x_n\arrow{d}{=}\\
    Mx_0\arrow{rr}{g}&&x_n
  \end{tikzcd}\]
  (i.e., the set of isomorphisms filling the rectangle), where $m$
  denotes the multiplication operation on $M$.
\item Composition is done in the obvious manner.
\end{itemize}
A monoid in this categorical theory is exactly an $M$-algebra in
$\cat{C}$.
\end{example}

\subsection{Theorization of theories}
\label{sec:theorize-theory}
\setcounter{subsubsection}{-1}

\subsubsection{}

In Section \ref{sec:theorize-category}, we have theorized the notion
of category by considering a category as a structure on the family of
sets consisting of the sets of maps.
Recall from Example \ref{ex:initially-graded} that a category was an
`initially graded' $1$-theory.
Example \ref{ex:graded-as-theorization} shows thus that category is a
theorization of $\Ini$-monoid, or set, in the sense
of Definition \ref{def:theorization-example}.
For these reasons, we shall call a categorical theory also an
``$\Ini$-graded $2$-theory''.

One can similarly consider the structure of a categorical theory as a
structure on the sets of its $2$-multimaps, and then try to theorize
the notion of categorical theory after fixing the collections of
objects and of $1$-morphisms.
It turns out that there is indeed an interesting theorization in this
case.
One might call the resulting theorized object a categorical $2$-theory
or an initially graded $3$-theory.

One might ask whether it is possible to iterate theorization in a
similar manner here, or starting not from set, but from some other
kind of structure.
We should of course ask possibility of non-trivial, and disirably,
interesting (in some way), theorizations.
In the case where the answer to the question is affirmative, just as
the original structure could be expressed as the structure of a monoid
over the terminal object among the theorized objects of the same kind,
the theorization similarly becomes the structure of a monoid over the
terminal object among the twice theorized objects, and so on, so all
the structures can be described, in a non-trivial manner, using their
iterated theorizations.
Moreover, by considering an algebra over a non-terminal theory (which
exists for a non-trivial theorization), an algebra over it, and so on,
one obtains various general structures, which can all be treated in a
unified manner.

The question asked above is non-trivial, since it asks the existence
of non-trivial theorizations.
See Section \ref{sec:theorization-general}.
However, we introduce the notion of $n$-theory in \cite{theory}, which
will inductively be an interesting theorization of ($n-1$)-theory.
The hierarcy as $n$ varies, of $n$-theories will be an infinite
hierarchy of iterated theorizations which extends the various standard
hierarchies of iterated categorifications, in particular, the hierarcy
of $n$-categories in the ``initially graded'' case.

While our higher theories will be in general a completely new
mathematical objects, we have already found very classical objects of
mathematics among $2$-theories.
Namely, while we have seen that a categorical theory was an
``initially graded'' $2$-theory, we have also noted in Section
\ref{sec:theorize-category}, that planar multicategories were
among categorical theories.

We have also seen a non-classical object among $2$-theories, as
Example \ref{ex:theory-from-monad}.
Less exotic examples of higher theories comes from the construction of
Section \ref{sec:basic-construction}.
Namely, a higher categorified instance of a \emph{lower} theorized
structure leads to a higher theorized structure through the iterated
application of the construction $\wasreteori$.
As an object, this is less interesting among the general higher
theories for the very reason that it is represented by a lower
theory.
However, the \emph{functors} between these theories are actually
interesting in that it is much more general than functors which we
consider between the original lower theories.
Namely, a functor between such higher theories amounts to highly
\emph{relaxed} functor between the higher categorified lower theories,
as follows from an iteration of the remark in Section
\ref{sec:basic-construction}.

In \cite{theory}, we discuss a general construction which we call
``delooping'' (see Section \ref{sec:further-development}), through
which we obtain an ($n+1$)-theory which normally fails to be
representable by an $n$-theory.
Another construction, which is closely related to the classical Day
convolution, will also be discussed in \cite{theory}, and this also
produces similar examples.

\subsubsection{}

The notion of $2$-theory immediately leads to a generalization of
Metaproposition \ref{metapro:enriching-multicategory}.
Indeed, for a symmetric multicategory $\cat{U}$, one can define a
theorization of $\cat{U}$-graded $1$-theory, which we call
\kore{$\cat{U}$-graded $2$-theory}.
It follows from the general idea on theorization that the notion of
$\cat{U}$-graded $1$-theory makes sense naturally in a
$\cat{U}$-graded $2$-theory, and this notion gives a generalization of
the ``metaproposition''.

Indeed, given a $\cat{U}\tensor E_1$-monoidal category $\cat{A}$, one
can categorically deloop $\cat{A}$ using the $E_1$-monoidal structure,
to obtain a $\cat{U}$-monoidal $2$-category $\deloop\cat{A}$ and
hence a $\cat{U}$-graded $1$-theory $\wasreteori\deloop\cat{A}$
enriched in categories.
Since this is a categorification of $\cat{U}$-graded $1$-theory (which
is only as much coloured as $\cat{U}$), the notion of $\cat{U}$-graded
operad in $\wasreteori\deloop\cat{A}$ makes sense, which naturally
generalizes from the case where $\cat{A}$ is symmetric monoidal, the
notion of $\cat{U}$-graded operad in $\cat{A}$.
Moreover, it is easy to check, when $\cat{U}$ is the few simplest
operads, that this coincides with the usual notion.
However, this notion of operad, including the coloured cases of it,
turns out to be nothing but the notion of $1$-theory in the
$\cat{U}$-graded $2$-theory $\wasreteori(\wasreteori\deloop\cat{A})$.

\subsubsection{}

More generally, with appropriate notion of grading, an $n$-theory
makes sense in an ($n+1$)-theory.
Let us briefly discuss the notion of grading for higher theories.

A starting point is that there is notion of algebra over each
$n$-theory.
We have mentioned that the monoid (i.e., algebra in the category of
sets) over the terminal $n$-theory coincides with the notion of
($n-1$)-theory.
For an $n$-theory $\cat{U}$, we can theorize the notion of
$\cat{U}$-monoid, generalizing from the case where $\cat{U}$ is
terminal.
We call our theorization \kore{$\cat{U}$-graded $n$-theory}, where our
term comes from the following.
Compare with our discussion in Section \ref{sec:theorize-algebra}.

\begin{proposition}
\label{prop:graded-is-overlying-intro}
A $\cat{U}$-graded $n$-theory is equivalent as data to a
\textbf{symmetric} $n$-theory $\cat{X}$ equipped with a functor
$\cat{X}\to\cat{U}$ of symmetric $n$-theories.
\end{proposition}

It would therefore seem natural to call a $\cat{U}$-monoid a
$\cat{U}$-graded \kore{($n-1$)-theory}, and indeed, there is a natural
notion of $\cat{U}$-graded \kore{$0$-theory} of which $\cat{U}$-monoid
is an ($n-1$)-th theorization.
It is therefore also natural to call the intermediate theorizations as
$\cat{U}$-graded \kore{$m$-theory} for $1\le m\le n-2$.

On the other hand, there is also a natural theorization of
$\cat{U}$-graded $n$-theory, which we of course call $\cat{U}$-graded
\kore{($n+1$)-theory}.
We obtain the following fundamental results.

\begin{theorem}\label{thm:graded-theory-is-monoid}
A $\cat{U}$-graded $n$-theory is equivalent as data to a monoid over
the ($n+1$)-theory $\wasreteori\cat{U}$.
A $\cat{U}$-graded ($n+1$)-theory is equivalent as
data to a $\wasreteori\cat{U}$-graded ($n+1$)-theory.
\end{theorem}

It follows that the natural notion of $\cat{U}$-graded $m$-theory for
$m\ge n+2$, is simply $\wasreteori^m_n\cat{U}$-graded $m$-theory,
where $\wasreteori^m_n$ denotes the ($m-n$)-fold iteration of the
construction $\wasreteori$.

We now obtain that $\cat{U}$-graded $m$-theory is a notion which
naturally makes sense in a $\cat{U}$-graded ($m+1$)-theory.
This is the enriched version of the notion.

\begin{example}
Let us denote the terminal $n$-theory by $\unity^n_\Com$, and let
$\diag\colon\cat{U}\to\unity^n_\Com$ be the unique functor.
Since every $n$-theory $\cat{V}$ is graded by $\unity^n_\Com$, one
obtains from this a $\cat{U}$-graded $n$-theory $\diag^*\cat{V}$.
A $\cat{U}$-algebra in $\diag^*\cat{V}$ is a (coloured) functor
$\cat{U}\to\cat{V}$ of $n$-theories.
\end{example}

\subsection{Further developments}
\label{sec:further-development}
\setcounter{subsubsection}{-1}

\subsubsection{}

Let us preview a few more highlights of our work.

\subsubsection{}

In addition to monoid over an $n$-theory, we also define the notion of
monoid over a monoid over an $n$-theory (if $n\ge 2$), monoid over
such a thing (if $n\ge 3$), and so on.
For example, for an $n$-theory $\cat{U}$, a $\cat{U}$-graded
($n-2$)-theory can be expressed as a monoid over the terminal
$\cat{U}$-monoid, among which the terminal one is such that a monoid
over it is exactly a $\cat{U}$-graded ($n-3$)-theory, and so on.

More generally, for a $\cat{U}$-graded $m$-theory $\cat{X}$, we define
the notion of $\cat{X}$-monoid, and more generally, of
$\cat{X}$-graded $\ell$-theory, as well as the notion of higher
theory graded over such a thing, and so on.
We can actually give very simple definitions of all these, using
Theorem \ref{thm:graded-theory-is-monoid} as a general principle.
We obtain a generalization of Proposition
\ref{prop:graded-is-overlying-intro} with these new notions.
We also obtain the enriched versions of all the notions.

\subsubsection{}

One can define the notion of $n$-theory enriched in a symmetric
monoidal category $\cat{A}$.
Indeed, one obtains a symmetric monoidal ($n+1$)-category
$\deloop^n\cat{A}$ by iterating the categorical delooping
construction, and since a symmetric monoidal ($n+1$)-category is an
($n+1$)-th categorification of commutative algebra, one can apply the
construction $\wasreteori$ $n+1$ times to $\deloop^n\cat{A}$, to
obtain an ($n+1$)-theory, which we shall denote by
$\wasreteori^{n+1}\deloop^n\cat{A}$.
Then an $n$-theory in this ($n+1$)-theory is a natural notion of
$n$-theory enriched in $\cat{A}$.
For example, the case $\cat{A}=\Set$ of this is the original notion of
$n$-theory.

We generalize the delooping construction used above for symmeric
monoidal higher category, to a certain construction $\ddeloop$ which
produces a symmetric $n$-theory from a symmetric ($n-1$)-theory.
This is a generalization of the categorical delooping in such a manner
that, for a symmetric monoidal category $\cat{A}$, there is a natural
equivalence
$\wasreteori^{n+1}\deloop^n\cat{A}=\ddeloop^n\wasreteori\cat{A}$,
which makes it natural to define for a symmetric multicategory
$\cat{U}$ which is not necessarily of the form $\wasreteori\cat{A}$,
an $n$-theory enriched in $\cat{U}$ as an $n$-theory in the
($n+1$)-theory $\ddeloop^n\cat{U}$.

Incidentally, if $\cat{U}$ is not of the form $\wasreteori\cat{A}$,
then $\ddeloop^n\cat{U}$ is usually \emph{not} representable by a
categorified $n$-theory.

\subsubsection{}

For a symmetric monoidal $n$-category $\cat{A}$, we construct a
certain symmetric monoidal ($n+1$)-category $A_n\cat{A}$ and a functor
$A_n\cat{A}\to\deloop^n\Set$ of symmetric monoidal ($n+1$)-categories,
which induces a functor
$\wasreteori^{n+1}A_n\cat{A}\to\wasreteori^{n+1}\deloop^n\Set$ of
($n+1$)-theories.
The use of this is the following.
Namely, while we have already mentioned that an $n$-theory $\cat{U}$
can be considered as an $n$-theory in
$\wasreteori^{n+1}\deloop^n\Set$, the construction above allows us
to understand a $\cat{U}$-graded $m$-theory as an appropriate lift of
the theory to $\wasreteori^{n+1}A_n\deloop^m\Set$.

\subsubsection{}

We also touch on more topics, such as the following.
\begin{itemize}
\item Pull-back and push-forward constructions which changes gradings,
  and their properties.
\item Some other basic constructions such as a construction for higher
  theories related to Day's convolution.
\item A few more general hierarchies of iterated theorizations,
  associated to systems of operations with multiple inputs \emph{and}
  multiple outputs, such as operations of `shapes' of bordisms as in
  various versions of a topological field theory.
\end{itemize}

The last topic leads to vast generalizations of the relevant versions
of a topological field theory.
It turns out that the generality of the notion allows very simple but
in a way exotic examples.

\subsubsection{}

We also enrich everything we consider in the work, in the Cartesian
symmetric monoidal infinity $1$-category of infinity groupoids,
instead of in sets.
Fortunately, this does not add any difficulty to the discussions.

\end{document}